\documentclass[12pt,fceqn]{article}
\usepackage{amsmath}
\usepackage{amsfonts}
\usepackage{cite}
\usepackage{amssymb,graphics,graphicx,subfigure,caption2,rotating}
\usepackage{amsmath, latexsym}
\usepackage[amsmath,thmmarks]{ntheorem}
\usepackage{color}
\numberwithin{equation}{section}
\newtheorem{theorem}{Theorem}[section]
\newtheorem{lemma}{Lemma}[section]
\newtheorem{remark}{Remark}[section]

\newtheorem{cor}{Corollary}[section]

\newcommand{\beq}{\begin{equation}}
\newcommand{\eeq}{\end{equation}}
\newcommand{\beqn}{\begin{eqnarray}}
\newcommand{\eeqn}{\end{eqnarray}}
\date{}
\begin{document}

\date{}
\title{The perturbation bound of the extended vertical linear complementarity problem\thanks{This research was supported by National Natural Science Foundation of China
(No.11961082,12071159).}}
\author{Shiliang Wu\thanks{Corresponding author: slwuynnu@126.com}, Wen
Li\thanks{liwen@scnu.edu.cn}, Hehui Wang\thanks{wanghehui1994@126.com}\\
{\small{\it $^{\dag}$School of Mathematics, Yunnan Normal University,}}\\
{\small{\it Kunming, Yunnan, 650500, P.R. China}}\\
{\small{\it $^{\ddag}$School of Mathematical Sciences, South China Normal University,}}\\
{\small{\it Guangzhou, Guangdong, 510631, P.R. China}}\\
{\small{\it $^{\S}$School of Mathematics and Statistics, Yunnan University,}}\\
{\small{\it Kunming, Yunnan, 650091, P.R. China}}\\
}
 \maketitle
\begin{abstract}
In this paper, we discuss the perturbation analysis of the extended
vertical linear complementarity problem (EVLCP). Under the assumption of the row
$\mathcal{W}$-property, several absolute and relative perturbation
bounds of EVLCP are given, which can be reduced to some existing results. Some numerical
examples are given to show the proposed bounds.


\end{abstract}

\section{Introduction}
Let $M_i\in \mathbb{R}^{n\times n}$ and $q_i\in
\mathbb{R}^{n},~i=0,1,...,k$. The following minimization equation is
to find $x\in \mathbb{R}^{n}$ such that
\begin{equation}\label{p1}
\min\{M_{0}x+q_{0}, M_{1}x+q_{1},\ldots, M_{k}x+q_{k}\}=0,
\end{equation}
where  $\min$ is the component minimum operator, which can be
reduced to the following two models :

\begin{equation}\label{eq:11}
\min\{x, Mx+q\}=0,
\end{equation}
where $M\in \mathbb{R}^{n\times n}$ and $q\in \mathbb{R}^{n}$, and

\begin{equation}\label{VLCP}
\min\{x, M_{1}x+q_{1},\ldots, M_{k}x+q_{k}\}=0.
\end{equation}

The model (\ref{p1}) is called the extended vertical linear
complementarity problem (EVLCP), denoted by EVLCP
($\mathbf{M},\mathbf{q}$) with $\mathbf{M}=(M_{0},M_{1},\ldots,
M_{k})$ and  $\mathbf{q}=(q_0,q_1,\cdots,q_k)$, the models
(\ref{eq:11}) and (\ref{VLCP}) are called the linear complementarity
problem (LCP) (e.g., see \cite{Cottle92, Murty88}) and the vertical
linear complementarity problem (VLCP) (e.g. see \cite{Cottle70,
Zhang09}), respectively.

A very important problem in the computational sciences is that how
the solution variation is when the data is perturbed. More
specifically, for the model (\ref{p1}), when $\Delta
M_{i}~\mbox{and}~\Delta q_i$ are the perturbation of
$M_{i}~\mbox{and}~q_i$, respectively, $i=0, 1,...,k$, how do we
characterize the change in the solution of the following perturbed
model:

\begin{equation}\label{p2}
 \min\{\tilde M_{0}y+\tilde q_{0}, \tilde M_{1}y+\tilde q_{1},\ldots, \tilde M_{k}y+\tilde q_{k}\}=0,
\end{equation}
where $\tilde M_{i}=M_{i}+\Delta M_{i},~\mbox{and}~\tilde
q_{i}=q_i+\Delta q_i, i=0,1,...,k.$ This problem has been
extensively studied for the model (\ref{eq:11}), which is often used
to deduce the sensitivity and stability analysis \cite{Cottle92} for
solving the LCP ($M,q$). A classical result on the error bound of
the LCP ($M,q$) was given in \cite{Mathias90} by Mathias and Pang.
In order to introduce the perturbation result for the model
(\ref{eq:11}), we first give some definitions and notations.

A matrix $M$ is called a $P$-matrix if all principal minors of $M$
are positive, in this case the model \eqref{eq:11} has the unique
solution (e.g., see  \cite{Cottle92}). Let
\[
c(M):=\min_{\|x\|_{\infty}=1}\big\{\max_{1\leq i\leq
n}x_{i}(Mx)_{i}\big\}.
\]
By the above definition, the perturbation bound for LCP (\ref{eq:11}) was given as follows.

\begin{lemma} [7.3.10 Lemma, \cite{Cottle92}]
For $M$ being a $P$-matrix, the following results hold:
\begin{description}
\item $(a)$ For any two vectors $q$ and $\check{q}$ in $\mathbb{R}^{n}$,
\[
\|x^{\ast}-x^{\star}\|_{\infty}\leq
\frac{1}{c(M)}\|q-\check{q}\|_{\infty},
\]
where $x^{\ast}$ and $x^{\star}$ denote the unique solutions of the
LCPs $(M, q)$ and $(M, \check{q})$, respectively.
\item $(b)$ For each vector $q\in \mathbb{R}^{n}$, there exist a neighborhood $U$ of the pair $(M,q)$ and a
constant $c_{0} > 0$ such that for any $(\bar{M}, \bar{q}),
(\hat{M}, \hat{q})\in U$, $\bar{M}, \hat{M}$ are $P$-matrices and
\[
\|x-y\|_{\infty}\leq
c_{0}(\|\bar{q}-\hat{q}\|_{\infty}+\|\bar{M}-\hat{M}\|_{\infty}),
\]
where $x$ and $y$ denote the unique solutions of the LCPs ($\bar{M},
\bar{q}$) and $(\hat{M}, \hat{q})$, respectively.
\end{description}
\end{lemma}

Lemma 1.1 is very important for studying the sensitivity and
stability analysis of the LCP model theoretically. Alternatively,
Chen and Xiang in \cite{Chen07} provided the following sharper
perturbation bounds than the ones in Lemma 1.1 by introduced the
following constant:
\[
\beta_{p}(M):=\max_{d\in [0,1]^{n}}\|(I-D+DM)^{-1}D\|_{p}
\]
for a $P$-matrix, where  $I$ is the identity matrix, $D$ is a
diagonal matrix whose diagonal entry is in $[0,1]$ and
$\|\cdot\|_{p}$ denotes the $p$-norm with $p\geq1$.
\begin{lemma} [Theorem 2.8 of \cite{Chen07}]
For $M$ being a $P$-matrix, the following results hold:
\begin{description}
\item $(a)$ For any two vectors $q$ and $\check{q}$ in $\mathbb{R}^{n}$,
\[
\|x^{\ast}-x^{\star}\|_{p}\leq \beta_{p}(M)\|q-\check{q}\|_{p},
\]
where $x^{\ast}$ and $x^{\star}$ denote the unique solutions of the
LCPs $(M, q)$ and $(M, \check{q})$, respectively.
\item $(b)$ For $\bar{M}, \hat{M}\in \mathcal{M}:=\{A~|~\beta_{p}(M)\|M-A\|_{p}\leq
\eta<1\}$ and $\bar{q},\hat{q}\in \mathbb{R}^{n}$,
\[
\|x-y\|_{p}\leq
\frac{\beta^{2}_{p}(M)}{(1-\eta)^{2}}\|(-\hat{q})_{+}\|_{p}\|\bar{M}-\hat{M}\|_{p}+\frac{
\beta_{p}(M)}{1-\eta}\|\bar{q}-\hat{q}\|_{p},
\]
where $x$ and $y$ denote the unique solutions of the LCPs ($\bar{M},
\bar{q}$) and $(\hat{M}, \hat{q})$, respectively.
\end{description}
\end{lemma}

It is difficult to compute $c(M)$, $c_{0}$ and $\beta_{p}(M)$. To
overcome this drawback, Chen and Xiang in \cite{Chen07} provided
some computable bounds for $M$ being an $H$-matrix with positive
diagonals, a symmetric positive definite matrix, a positive definite
matrix, respectively.

The algorithm, applications and the existence of solutions for the
model (\ref{p1}) have been given (see, e.g., \cite{Fujisawa72,
Sun89, Oh86, Zabaljauregui21, Goeleven96, Gowda96, Zhou08,
Habetler92,Gowda, Qi99, Gowda962, Ebiefung93, Mezzadri}). However,
so far, to our knowledge, the perturbation analysis of models
(\ref{p1}) and (\ref{VLCP}) has not been discussed. In order to fill
in this study gap, in this paper, inspired by the success as in LCP
model \eqref{eq:11} as in \cite{Chen07}, we focus on discussing the
perturbation analysis for EVLCP (\ref{p1}) and VLCP (\ref{VLCP}).
The contributions are given below:
\begin{itemize}
\item The framework of perturbation bound for the EVLCP model is proposed by developing the
technique given in \cite{Chen07}, from which one may derive the
corresponding bound given in \cite{Chen07}.
\item Some computable bounds are also given for some special block matrices $\mathbf{M}=(M_{0},M_{1},$ $\ldots, M_{k})$.
\item Some examples from discretization of Hamilton-Jacobi-Bellman (HJB) equation are given to show the proposed
bounds.
\end{itemize}

The rest of the article is organized as follows. Section 2 is
preliminary, in which we give some lemmas. In Section 3, we provide
some perturbation bounds for the EVLCP ($\mathbf{M},\mathbf{q}$) by
using a general equivalent form of the minimum function under the
row $\mathcal{W}$-property. In Section 4, some relative perturbation
bounds are provided. In Section 5, some numerical examples are given
to show the feasibility of the perturbation bound. Finally, in
Section 6, we give some conclusion remarks to end this paper.

\vspace{0.3cm}
Finally, in this section we give some notations and definitions \cite{Berman, Gowda}, which will be used in the sequel.

Let $A = (a_{ij})$, $B = (b_{ij}) \in \mathbb{R}^{n\times n}$ and
$N=\{1,2,\ldots,n\}$. Then we denote $|A|=(|a_{ij}|)$. The order $A\geq
(>)B$ means $a_{ij}\geq (>) b_{ij}$ for any $i,j\in N$.

$A =(a_{ij})$ is called an $M$-matrix if $A^{-1}\geq0$ and $a_{ij}\leq0$
($i\neq j$) for $i,j \in N$; an $H$-matrix if its comparison matrix $\langle A\rangle$ (i.e., $\langle a\rangle_{ii}=|a_{ii}|, \langle
a\rangle_{ij}=-|a_{ij}|$  $i\neq j$ for $i,j \in N$) is an $M$-matrix; an $H_{+}$-matrix if $A$ is an $H$-matrix
with $a_{ii} > 0$ for  $i\in N$; a strictly diagonal dominant (sdd) matrix if $|a_{ii}|>\sum_{j\neq i}|a_{ij}|,\ i\in N$;
an irreducible diagonal dominant (idd) matrix if $A$ is an irreducible, $|a_{ii}|\geq\sum_{j\neq i}|a_{ij}|,\ i\in N$ and $\{i\in N: |a_{ii}|>\sum_{j\neq i}|a_{ij}|\}\neq\varnothing$.

Let $e=(1,1,\ldots,1)^{T}$ and by $\rho(\cdot)$ we denote the
spectral radius of a matrix. For a vector $q\in \mathbb{R}^{n}$, by
$q_{+}$ and $q_{-}$ we denote $q_{+}=\max\{0,q\}$ and
$q_{-}=\max\{0,-q\}$.

In this paper, the norm $\|\cdot\|_{p}$ means $p$-norm with
$p\geq1$.

Let
\[\mathcal{D}=\{(D_{0},D_{1},\ldots,D_{k})~|~D_{i}=\mbox{diag}(d_{i})\
\mbox{with} \ d_{i}\in [0,1]^{n} \ (i=0,1,\ldots,k)\ \mbox{and}\
\sum_{i=0}^{k}D_{i}=I\}.\]

A block matrix $\mathbf{M}=(M_{0},M_{1},\ldots, M_{k})$ is said to be with {\bf row $\mathcal{W}$-property} if
\[
\min(M_{0}x, M_{1}x,\ldots, M_{k}x)\leq0\leq \max(M_{0}x,
M_{1}x,\ldots, M_{k}x)\Rightarrow x=0.
\]

It is noted that EVLCP ($\mathbf{M},\mathbf{q}$) has the unique
solution for any $\mathbf{q}$ if and only if $\mathbf{M}$ has the
row $\mathcal{W}$-property (see Theorem 17 in \cite{Gowda}).

\section{Some lemmas}
In this section, we give some basic lemmas, which will be used in
the sequel. The first one is a general equivalent formula of the
minimum function.
\begin{lemma}
Let all $a_{i}, b_{i}\in \mathbb{R}$, $i=1,2,\ldots,n$. Then there
exist $\lambda_{i}\in [0,1]$ with $\sum_{i=1}^{n}\lambda_{i}=1$ such
that
\begin{equation}\label{eq:21}
\min_{1\leq i\leq n}\{a_{i}\}-\min_{1\leq i\leq
n}\{b_{i}\}=\sum_{i=1}^{n}\lambda_{i}(a_{i}-b_{i}).
\end{equation}
\end{lemma}
\textbf{Proof.} The result follows immediately from the mean value
theorem of Lipschitz functions with the generalized gradient (see
2.3.7 Theorem in \cite{Clarke}). $\hfill{} \Box$

\begin{lemma}
The block matrix $\mathbf{M}=(M_{0}, M_{1},\ldots, M_{k})$ has the row
$\mathcal{W}$-property if and only if
$D_{0}M_{0}+D_{1}M_{1}+\ldots+D_{k}M_{k}$ is nonsingular for any
$(D_{0},D_{1},\ldots,D_{k})\in \mathcal{D}$.
\end{lemma}
\textbf{Proof.} The lemma follows immediately from Lemma 2.1 and Theorem 3 (b) in
\cite{Sznajder95}. $\hfill{} \Box$

\begin{lemma}
Let $a_{i}, b_{i}, t_{i}\in \mathbb{R}$ with $a_{i}>0$, $t_{i}\in
[0,1]$ and $\sum_{i=1}^{n}t_{i}=1$. Then
\begin{equation*}
\frac{\sum_{i=1}^{n}t_{i}b_{i}}{\sum_{i=1}^{n}t_{i}a_{i}}\leq\max\bigg\{\frac{|b_{i}|}{a_{i}}\bigg\}.
\end{equation*}
\end{lemma}
 {\textbf{Proof.} The result follows immediately from the Cauchy inequality. $\hfill{} \Box$

\begin{lemma}
Let $a,b\geq 0$. Then
\[
\frac{t}{ta+(1-t)b}\leq\frac{1}{a}
\]
for any $t\in[0,1]$. 
\end{lemma}
 \textbf{Proof.} Let
\[
f(t)=\frac{t}{ta+(1-t)b}.
\]
Then
\[
f'(t)=\frac{b}{(ta+(1-t)b)^{2}}>0,
\]
it follows that $f(t)$ is a strictly monotone increasing function of
$t$ with $t\in[0,1]$. Therefore,  when $t=1$, we obtain that
$f(t)_{\max}=\frac{1}{a}$. $\hfill{} \Box$

\section{Perturbation bounds}

In this section, we always assume that the block matrix $\mathbf{M}$
has the row $\mathcal{W}$-property without further illustration. In
this case, EVLCP ($\mathbf{M},\mathbf{q}$) has the unique solution
$x^*$.
\subsection{Framework of perturbation analysis for EVLCP}
In this subsection, we discuss perturbation analysis of EVLCP when
both $\mathbf{M}$ and $\mathbf{q}$ are perturbed to
$\tilde{\mathbf{M}}=(\tilde{M}_0,\tilde{M}_1,\cdots,\tilde{M}_k)$
and
$\tilde{\mathbf{q}}=(\tilde{q}_0,\tilde{q}_1,\cdots,\tilde{q}_k)$,
respectively. Assume that $\tilde{\mathbf{M}}$ has the row
$\mathcal{W}$-property. Then EVLCP
($\tilde{\mathbf{M}},\tilde{\mathbf{q}}$) has the unique solution
$y^*$.

Let $r_{i}=M_ix^*+q_i:=(a_1^{(i)},...,a_n^{(i)})^T$ and
$\tilde{r}_{i}=
\tilde{M}_iy^*+\tilde{q}_i:=(b_1^{(i)},...,b_n^{(i)})^T$. It follows
from Lemma 2.1 that there exist $d_{\ell}^{(i)}\in [0,1]$ with
$\sum_{i=0}^{k} d_{\ell}^{(i)}=1$ ($\ell=1,2,\ldots,n$) such that
$$
\min_{0\leq i\leq k} \{a_{\ell}^{(i)}\}-\min_{0\leq i\leq
k}\{b_{\ell}^{(i)}\}=\sum_{i=0}^{k}d_{\ell}^{(i)}(a_{\ell}^{(i)}-b_{\ell}^{(i)}).
$$
Let $\tilde
D_{i}=\mbox{diag}(d_{1}^{(i)},...,d_{n}^{(i)}),~i=0,1,...,k$. Then
it is easy to see that for each $i$,  $\tilde{D}_{i}$ is a
nonnegative diagonal matrix, and $\sum_{i=0}^{k}\tilde{D}_{i}=I$
such that
\begin{equation}\label{p11}
  \min\{r_0,r_1,...,r_k\}-\min\{\tilde{r}_{0},\tilde{r}_{1},...,\tilde{r}_{k}\}=\sum_{i=0}^{k}\tilde{D}_i(r_{i}-\tilde
r_{i}).
\end{equation}

Since $x^*$ and $y^*$ are the solution of EVLCP
($\mathbf{M},\mathbf{q}$) and EVLCP
($\tilde{\mathbf{M}},\tilde{\mathbf{q}}$), respectively, we have
$$\sum_{i=0}^{k}\tilde{D}_i(r_{i}-\tilde r_{i})=0,$$ which implies
that
\begin{eqnarray*}
(\sum_{i=0}^{k}\tilde{D}_iM_i)x^*=(\sum_{i=0}^{k}\tilde{D}_i\tilde{M}_i)y^*+\sum_{i=0}^{k}\tilde{D}_i(\tilde{q}_i-q_i).
\end{eqnarray*}
Hence we get
\begin{eqnarray}\label{eq:31}
(\sum_{i=0}^{k}\tilde{D}_iM_i)(x^*-y^*)=(\sum_{i=0}^{k}\tilde{D}_i(\tilde{M}_i-M_i))y^*+\sum_{i=0}^{k}\tilde{D}_i(\tilde{q}_{i}-q_i).
\end{eqnarray}
Let $\tilde{S}_{\textbf{M}}=\sum_{i=0}^{k}\tilde{D}_iM_i$. It follows from Lemma 2.2 that $\tilde{S}_{\textbf{M}}$ is nonsingular, and then by Eq. (\ref{eq:31}), we have
\begin{align}\label{p2}
 x^*-y^*=&\tilde{S}_{\textbf{M}}^{-1}\big[(\sum_{i=0}^{k}\tilde{D}_i(\tilde{M}_i-M_i))y^*+\sum_{i=0}^{k}\tilde{D}_i(\tilde{q}_{i}-q_i)\big] \nonumber \\
=&\sum_{i=0}^{k}\tilde{S}_{\textbf{M}}^{-1}\tilde{D}_i[(\tilde{M}_i-M_i)y^*+(\tilde{q}_{i}-q_i)].
\end{align}
For any $(D_{0},D_{1},\ldots,D_{k})\in \mathcal{D}$, setting
$$S_{\textbf{M}}=\sum_{i=0}^{k}D_iM_i$$ and
\begin{align} \label{eq:32}
\alpha_{i}(\textbf{M})=\max_{(D_{0},D_{1},\ldots,D_{k})\in \mathcal{D}}\|S_{\textbf{M}}^{-1}D_i\|, i=0,1,\ldots,k.
\end{align}
Then
\begin{align*}
\alpha_{i}(\textbf{M})\geq
\|\tilde{S}_{\textbf{M}}^{-1}\tilde{D}_i\|, i=0,1,\ldots,k.
\end{align*}
So, by (\ref{p2}) and (\ref{eq:32}) we have
\begin{align}\label{eq:33}
\|x^*-y^*\|
\leq\sum_{i=0}^{k}\alpha_{i}(\textbf{M})\|(\tilde{M}_i-M_i)y^*+(\tilde{q}_{i}-q_i)\|,
\end{align}
which gives a perturbation bound for EVLCP.

\begin{theorem} Let $\mathbf{M}=(M_{0},M_{1},\ldots, M_{k})$
have the row $\mathcal{W}$-property and $\alpha_{i}(\emph{\textbf{M}})$ be defined as in $(\ref{eq:32})$.
Then for both block vectors $\mathbf{q}=(q_0,q_1,\cdots,q_k)$ and
$\tilde{\mathbf{q}}=(\tilde{q}_0,\tilde{q}_1,\cdots,\tilde{q}_k)$,
EVLCP $(\mathbf{M},\mathbf{q})$ and EVLCP
$(\mathbf{M},\tilde{\mathbf{q}})$ have the unique solution $x^*$ and
$y^*$, respectively. Furthermore,

\begin{eqnarray}\label{eq:34}
\|x^*-y^*\|\leq
\sum_{i=0}^{k}\alpha_{i}(\emph{\textbf{M}})\|\tilde{q}_{i}-q_i\|.
\end{eqnarray}

\end{theorem}
\begin{remark}
Taking $k=1~ \mbox{and} ~M_0=I$, it is easy to see that
$\beta_p(M)=\alpha_1(\emph{\textbf{M}})$. Furthermore, taking
$q_0=0, \tilde{q_0}=0$, the bound \eqref{eq:34} reduces to the
corresponding one in \cite{Chen07}.
\end{remark}

Let $\mathbf{M}=(M_{0},M_{1},\ldots, M_{k})$ have the row $\mathcal{W}$-property.
Then we define $\mathfrak{M}_\eta$ as follows:
$$\mathfrak{M}_\eta:=\{\mathbf{A}=(A_0,A_1,\cdots,A_k)~|~\sum_{i=0}^{k}\alpha_{i}(\textbf{M})\|A_i-M_i\|\leq
\eta< 1\},$$ where $\alpha_{i}(\textbf{M})$ is defined as in
$(\ref{eq:32})$.

In order to get a general perturbation bound, it needs the following
lemma:
\begin{lemma}\label{3.1} The following statements hold:
\begin{description}
\item  $ (1)$ A block  matrix $\mathbf{A}=(A_0,A_1,\cdots,A_k)\in
\mathfrak{M}_\eta$ is with the row $\mathcal{W}$-property.
\item  $ (2)$  For any $\mathbf{A}\in \mathfrak{M}_\eta$,
$\alpha_{i}(\textbf{\emph{A}})\leq\tau_{i}(\textbf{\emph{M}})$,
where
\[
\tau_{i}(\textbf{\emph{M}})=\frac{1}{1-\eta}\alpha_{i}(\textbf{\emph{M}}).
\]
\end{description}
\end{lemma}
\textbf{Proof.} First, we show the assertion (1) holds. For any
$\mathbf{A}=(A_0,A_1,\cdots,A_k)\in \mathfrak{M}_\eta$, since
\begin{align*}
  \|S_{\textbf{M}}^{-1}\sum_{i=0}^{k}D_i(A_i-M_i)\| =& \|S_{\textbf{M}}^{-1}D_0(A_0-M_0)+S_{\textbf{M}}^{-1}D_1(A_1-M_1)\\
  &+\ldots+S_{\textbf{M}}^{-1}D_k(A_k-M_k)\|\\
  \leq&  \|S_{\textbf{M}}^{-1}D_0\|\|A_0-M_0\|+\|S_{\textbf{M}}^{-1}D_1\|\|A_1-M_1\|\\
  &+\ldots+\|S_{\textbf{M}}^{-1}D_k\|\|A_k-M_k\|\\
  \leq& \sum_{i=0}^{k}\alpha_{i}(\textbf{M})\|A_i-M_i\|\leq \eta< 1
\end{align*}
and
\begin{equation}\label{eq:35}
 \sum_{i=0}^{k}D_iA_i=S_{\textbf{M}}[I+S_{\textbf{M}}^{-1}\sum_{i=0}^{k}D_i(A_i-M_i)],
\end{equation}
it is known that the matrix $\sum_{i=0}^{k}D_iA_i$ is nonsingular for any $(D_{0},D_{1},\ldots,D_{k})\in \mathcal{D}$.
By Lemma 2.2, $\mathbf{A}$ has the row $\mathcal{W}$-property. This proves the assertion (1).

Next, we will prove the assertion (2). From (\ref{eq:35}) we have
\begin{align*}
 (\sum_{i=0}^{k}D_iA_i)^{-1}D_i=[I+S_{\textbf{M}}^{-1}\sum_{i=0}^{k}D_i(A_i-M_i)]^{-1}S_{\textbf{M}}^{-1}D_i,
\end{align*}
which together with the definition of $S_{\textbf{A}}$ gives
\begin{align*}
 \|S_{\textbf{A}}^{-1}D_i\|\leq\|[I+S_{\textbf{M}}^{-1}\sum_{i=0}^{k}D_i(A_i-M_i)]^{-1}\|\|S_{\textbf{M}}^{-1}D_i\|.
\end{align*}
It is easy to see that
\begin{align*}
\|[I+S_{\textbf{M}}^{-1}\sum_{i=0}^{k}D_i(A_i-M_i)]^{-1}\|\leq\frac{1}{1-\sum_{i=0}^{k}\alpha_{i}(\textbf{M})\|A_i-M_i\|}\leq
\frac{1}{1-\eta}.
\end{align*}
Hence we have
\begin{align}\label{36}
 \|S_{\textbf{A}}^{-1}D_i\|\leq \frac{1}{1-\eta}\|S_{\textbf{M}}^{-1}D_i\|.
\end{align}
The desired bound follows from \eqref{36}. This proves (2).
$\hfill{} \Box$

\begin{lemma}\label{3.2} Let $\mathbf{M}=(M_{0},M_{1},\ldots, M_{k})$
have the row $\mathcal{W}$-property, and
$\mathbf{q}=(q_0,q_1,\cdots,q_k)$. Then there exists $(\tilde
D_{0},\tilde D_{1},\ldots,\tilde D_{k})\in \mathcal{D}$ such that
the unique solution $x^*$ of EVLCP $(\mathbf{M},\mathbf{q})$ is
given by
\begin{eqnarray}\label{37}
x^*=-\sum_{i=0}^{k}\tilde S_{\textbf{M}}^{-1}\tilde D_{i}q_i,
\end{eqnarray}
where $\tilde{S}_{\textbf{M}}=\sum_{i=0}^{k}\tilde{D}_iM_{i}.$
\end{lemma}
\textbf{Proof.} Let
$r_{i}=M_ix^*+q_i:=(a_1^{(i)},...,a_n^{(i)})^T,~i=0,1,...,k$. Since
$x^*$ is the solution of the EVLCP ($\mathbf{M},\mathbf{q}$), we
have $\min\{r_{0},r_{1},\ldots, r_{k}\}=0.$ This implies that for
any $s$ there is a vector $r_{i}$ whose the $s$-th component is
equal to zero, $s=0,1,...,k$. Assume that the number of vectors in
$\{r_{0}, r_{1},...,r_{k}\}$ whose the $s$-th component is $0$ is
$t_s$, say $a_s^{(i_1)}=a_s^{(i_2)}=...=a_s^{(i_{t_s})}=0$. Taking
non-negative diagonal matrices
$\tilde{D}_i=\mbox{diag}(d_1^{(i)},...,d_n^{(i)})$ such that
$d_{s}^{(i_1)}=...=d_{s}^{(i_{t_s})}=\frac{1}{t_s}$ and
$d_{s}^{(i)}=0,~i\neq i_1,...,i_{t_s}$. Then it is easy to see that
$(\tilde D_{0},\tilde D_{1},\ldots,\tilde D_{k})\in \mathcal{D}$ and
\begin{align*}
 \sum_{i=0}^{k}(\tilde{D}_i (M_ix^*+q_i))=\sum_{i=0}^{k}\tilde{D}_ir_i=0.
\end{align*}
This implies
\begin{align*}
 (\sum_{i=0}^{k}\tilde{D}_i M_{i})x^*=-\sum_{i=0}^{k}\tilde{D}_iq_{i},
\end{align*}
from which it follows that
\begin{align*}
 x^*=-(\sum_{i=0}^{k}\tilde{D}_iM_{i})^{-1}\sum_{i=0}^{k}\tilde{D}_iq_{i}.
\end{align*}
This proves the lemma. $\hfill{} \Box$

The following theorem is the framework of EVLCP perturbation.

\begin{theorem}
Let both $\mathbf{A}=(A_0,A_1,\cdots,A_k)$ and
$\mathbf{B}=(B_0,B_1,\cdots,B_k)$ be in $\mathfrak{M}_\eta$, and let
$\mathbf{\bar{q}}=(\bar{q}_0,\bar{q}_1,\cdots,\bar{q}_k)$ and
$\mathbf{\bar{p}}=(\bar{p}_0,\bar{p}_1,\cdots,\bar{p}_k)$. Then
EVLCP $(\mathbf{A},\mathbf{\bar{q}})$ and EVLCP
$(\mathbf{B},\mathbf{\bar{p}})$ have the unique solutions $x^*$ and
$y^*$, respectively. Moreover, we have
\begin{eqnarray}\label{eq:38}
\|x^*-y^*\|\leq
\bigg(\sum_{i=0}^{k}\tau_{i}(\textbf{\emph{M}})\|A_i-B_i\|\bigg)\bigg(\sum_{i=0}^{k}\tau_{i}(\textbf{\emph{M}})\|\bar{p}_i\|\bigg)+\sum_{i=0}^{k}\tau_{i}\|\bar{q}_i-\bar{p}_i\|,
\end{eqnarray}
where $\tau_{i}(\textbf{\emph{M}})$ is given by Lemma \ref{3.1}.
\end{theorem}
\textbf{Proof.} The first assertion follows from Lemma 3.1. Next we show the perturbation bound. 
By (\ref{eq:33}) we have
\begin{eqnarray}\label{eq:39}
\|x^*-y^*\|\leq
\sum_{i=0}^{k}\alpha_{i}(\textbf{A})(\|(B_i-A_i)\|\|y^*\|+\|\bar{p}_i-\bar{q}_i\|).
\end{eqnarray}
It follows from Lemma 3.2 that
\begin{align*}
 y^*=-\sum_{i=0}^{k}\tilde{S}^{-1}_{\textbf{B}}\tilde{D}_i\bar{p}_{i},
\end{align*}
where $\tilde{S}_{\textbf{B}}=\sum_{i=0}^{k}\tilde{D}_iB_{i}$.
Clearly we have
\begin{align*}
\alpha_{i}(\textbf{B})\geq
\|\tilde{S}_{\textbf{B}}^{-1}\tilde{D}_i\|, i=0,1,\ldots,k.
\end{align*}
Hence,
\begin{align}\label{eq:310}
\|y^*\|\leq\sum_{i=0}^{k}\alpha_{i}(\textbf{B})\|\bar{p}_{i}\|.
\end{align}
By (\ref{eq:39}), (\ref{eq:310}) and Lemma 3.1, we obtain the
desired bound (\ref{eq:38}).  $\hfill{} \Box$

\vspace{0.3cm} Next, we discuss perturbation analysis of VLCP
($\mathbf{M},\mathbf{q}$), i.e., $M_{0}=I$ and $q_{0}=0$ for Eq.
(\ref{p1}).

Let
\[
\tilde{S}^{v}_{\textbf{M}}=\tilde{D}_{0}+\sum_{i=1}^{k}\tilde{D}_iM_i,
S^{v}_{\textbf{M}}=D_{0}+\sum_{i=1}^{k}D_iM_i
\]
and
\[
\tilde{S}^{v}_{\textbf{B}}=\tilde{D}_{0}+\sum_{i=1}^{k}\tilde{D}_iB_i,
S^{v}_{\textbf{B}}=D_{0}+\sum_{i=1}^{k}D_iB_i,
\]
where matrices $\tilde{D}_{i}=\mbox{diag}(\tilde{d}_{i})$ with
$\tilde{d}_{i}\in [0,1]^{n}$ $(i=0,1,\ldots,k)$ are nonnegative
diagonal and $\sum_{i=0}^{k}\tilde{D}_{i}=I$, matrices
$D_{i}=\mbox{diag}(d_{i})$ with $d_{i}\in [0,1]^{n}$
$(i=0,1,\ldots,k)$ are arbitrary nonnegative diagonal and
$\sum_{i=0}^{k}D_{i}=I$. Then we take
\begin{align}\label{eq:311}
\tilde{\alpha}_{i}(\textbf{M})=\max_{(D_0, D_{1},\ldots, D_{k})\in
\mathcal{D}}\|(S^{v}_{\textbf{M}})^{-1}D_i\|,
\tilde{\alpha}_{i}(\textbf{B})=\max_{(D_0, D_{1},\ldots, D_{k})\in
\mathcal{D}}\|(S^{v}_{\textbf{B}})^{-1}D_i\|, i=1,\ldots,k.
\end{align}
Clearly,
\begin{align*}
\tilde{\alpha}_{i}(\textbf{M})\geq\|(\tilde{S}^{v}_{\textbf{M}})^{-1}\tilde{D}_i\|,
\tilde{\alpha}_{i}(\textbf{B})
\geq\|(\tilde{S}^{v}_{\textbf{B}})^{-1}\tilde{D}_i\|, i=1,\ldots,k.
\end{align*}

Notice that $0$ is the solution of the VLCP
($\mathbf{B},\mathbf{p_+}$) and $y^*$ is the unique solution of the
VLCP ($\mathbf{B},\mathbf{p}$). Then we have
\begin{align*}
\|y^*-0\|\leq\sum_{i=0}^{k}\tilde{\alpha}_{i}(\textbf{B})\|(-p_i)_+\|.
\end{align*}

By analogical proof to Theorem 3.1 and Theorem 3.2, we can obtain
the following theorem.

\begin{theorem} Let $\mathbf{M}=(I,M_{1},\ldots, M_{k})$
have the row $\mathcal{W}$-property and
$\tilde{\alpha}_{i}(\emph{\textbf{M}})$ be defined as in
$(\ref{eq:311})$. Then the following statements hold:
\begin{description}
\item  $ (i)$ For any two block vectors
$\mathbf{q}=(0,q_1,\cdots,q_k)$,
$\tilde{\mathbf{q}}=(0,\tilde{q}_1,\cdots,\tilde{q}_k)$ with
$q_i,\tilde{q}_i\in \mathbb{R}^n$,
\begin{eqnarray*}
\|x^*-y^*\|\leq
\sum_{i=1}^{k}\tilde{\alpha}_{i}(\textbf{\emph{M}})\|\tilde{q}_{i}-q_i\|,
\end{eqnarray*}
where $x^*$ and $y^*$ are the unique solutions of VLCP
$(\mathbf{M},\mathbf{q})$ and VLCP
$(\mathbf{M},\tilde{\mathbf{q}})$, respectively.
\item  $ (ii)$ Every block  matrix $\tilde{\mathbf{A}}=(I,A_1,\cdots,A_k)\in \tilde{\mathfrak{M}}:=\{\tilde{\mathbf{A}}|\sum_{i=1}^{k}\tilde{\alpha}_{i}(\textbf{\emph{M}})\|A_i-M_i\|\leq \eta< 1\}$
has the row $\mathcal{W}$-property. Let
\[
\tilde{\tau}_{i}(\textbf{\emph{M}})=\frac{1}{1-\eta}\tilde{\alpha}_{i}(\textbf{\emph{M}}).
\]
Then for any
$\tilde{\mathbf{A}}=(I,A_1,\cdots,A_k),\tilde{\mathbf{B}}=(I,B_1,\cdots,B_k)\in
\tilde{\mathfrak{M}}$,
$\mathbf{\hat{q}}=(0,\hat{q}_1,\cdots,\hat{q}_k)$ and
$\mathbf{\hat{p}}=(0, \hat{p}_1,\cdots,\hat{p}_k)$,
\begin{eqnarray*}
\|x^*-y^*\|\leq
\bigg(\sum_{i=1}^{k}\tilde{\tau}_{i}(\textbf{\emph{M}})\|A_i-B_i\|\bigg)\bigg(\sum_{i=1}^{k}\tilde{\tau}_{i}(\textbf{\emph{M}})\|(-\hat{p}_i)_{+}\|\bigg)+\sum_{i=1}^{k}\tilde{\tau}_{i}(\textbf{\emph{M}})\|\hat{q}_i-\hat{p}_i\|,
\end{eqnarray*}
where $x^*$ and $y^*$ are the unique solutions of the VLCP
$(\tilde{\mathbf{A}},\mathbf{\hat{q}})$ and the VLCP
$(\tilde{\mathbf{B}},\mathbf{\hat{p}})$, respectively.
\end{description}
\end{theorem}

\begin{remark}
Here we consider a special case. It is noted that $\mathbf{M}=(I,M)$
has the row $\mathcal{W}$-property if and only if $M$ is a
$P$-matrix (see page 696 of \cite{Sznajder95}). By Remark 3.1,
$\beta_p(M)=\alpha_1(\textbf{M})$,  the bounds in Theorem 3.3
reduces to the corresponding ones in \cite{Chen07} (also see Lemma
1.2).
\end{remark}

\subsection{Estimations of $\alpha_{i}(M)$ and $\tilde{\alpha}_{i}(M)$}
Under assumption of the row $\mathcal{W}$-property,  we have given
perturbation bounds of EVLCP and VLCP  in Theorems 3.2 and 3.3,
respectively. However,  it is difficult to compute quantities
$\alpha_{i}(\textbf{M})$ and $\tilde{\alpha}_{i}(\textbf{M})$. In
this subsection, we explore computability estimations for
$\alpha_{i}(\textbf{M})$ and $\tilde{\alpha}_{i}(\textbf{M})$ for
the special block matrix $\mathbf{M}$.

To calculate $\alpha_{i}(\textbf{M})$ and
$\tilde{\alpha}_{i}(\textbf{M})$, we consider two types of special
matrices: (1) all the diagonal entries of the matrices $M_{i}$ in
$\mathbf{M}$ are positive; (2) all the matrices $M_{i}$ in
$\mathbf{M}$ are an sdd  matrix.

\subsubsection{Case 1}
Here we consider the case (1), i.e., all the diagonal entries of
matrices $M_{i}$ in $\mathbf{M}$ are positive. Let $A_i\in
\mathbb{R}^{n\times n},~i=0,1,...,k$. We denote by $\max_{0\leq
i\leq k}\{A_{i}\}$ a matrix whose the $(i,j)$-entry is the largest
$(i,j)$-entry among all matrices $A_i\in \mathbb{R}^{n\times
n},~i=0,1,...,k$.

Let $\wedge_{i}$ be the diagonal part of $M_{i}$, and let
$M_{i}=\wedge_{i}-C_{i}$, $i=0,1,\ldots,k$. First we give a lemma.
\begin{lemma}
Let  all the diagonal elements of the matrices $M_{i}$ in $\mathbf{M}$ be positive. If
\begin{equation}\label{eq:312}
\rho(\max_{0\leq i\leq k}\{\wedge^{-1}_{i}|C_{i}|\})<1.
\end{equation}
Then 
$\mathbf{M}=(M_{0},M_{1},\ldots, M_{k})$ has the row $\mathcal{W}$-property.
\end{lemma}
\textbf{Proof.}  Let $V=\sum_{i=0}^{k}D_{i}\wedge_{i}$ and
$U=\sum_{i=0}^{k}D_{i}C_{i}$. Then it is easy to see that the matrix
$V$ is nonsingular, and thus
\[
S_{\textbf{M}}=V-U=V(I-V^{-1}U).
\]
From Lemma 2.3, we have
\[
V^{-1}U\leq V^{-1}|U|\leq \max_{0\leq i\leq
k}\{\wedge^{-1}_{i}|C_{i}|\}.
\]
By Theorem 8.1.18 of \cite{Horn}, we get
\[
\rho(V^{-1}U)\leq \rho(V^{-1}|U|)\leq \rho(\max_{0\leq i\leq
k}\{\wedge^{-1}_{i}|C_{i}|)\},
\]
which implies that $S_{\textbf{M}}$ is nonsingular. Hence, $\mathbf{M}=(M_{0},M_{1},\ldots, M_{k})$ has the row $\mathcal{W}$-property.  $\hfill{} \Box$

\begin{remark} It is known that all the diagonal entries are positive for some special matrices, e.g., each $M_{i}$ in
$\mathbf{M}$ is a symmetric positive definite matrix, an $M$-matrix
or an $H_{+}$-matrix. However the condition (\ref{eq:312}) in Lemma
3.3 can not be omitted,  a counter-example is given below. Taking
$k=1$, and
\end{remark}

\[M_{0}=\left[\begin{array}{cc}
2&1\\
1&1
\end{array}\right], M_{1}=\left[\begin{array}{cc}
1&1\\
1&2
\end{array}\right].
\]
Obviously, both $M_{0}$ and  $M_{1}$ are symmetric positive definite
matrices, and also are  $H_{+}$-matrices or idd matrices. If we take
\[\hat{D}_{0}=\left[\begin{array}{cc}
0&0\\
0&1
\end{array}\right],\hat{D}_{1}=\left[\begin{array}{cc}
1&0\\
0&0
\end{array}\right],
\]
then the matrix
\[\hat{D}_{0}M_{0}+\hat{D}_{1}M_{1}=\left[\begin{array}{cc}
1&1\\
1&1
\end{array}\right]
\]
is singular.

\begin{remark}  Remark 3.3 implies that the results in Theorem 2.5 and
Theorem 2.7 in \cite{Chen07} do not extend to EVLCP
$(\mathbf{M},\mathbf{q})$. That is to say, the general forms of
Theorem 2.5 and Theorem 2.7 in \cite{Chen07} are no longer valid,
e.g.,  taking $M_{0}$ and $M_{1}$ in Remark 3.3,
\[
\beta_{p}(M)=\max_{d\in [0,1]^{n}}\|((I-D)M_{0}+DM_{1})^{-1}D\|_{p}
\]
does not exist.
\end{remark}

Next theorem shows the computability of $\alpha_{\ell}(\emph{\textbf{M}})$.

\begin{theorem}
Under the assumption of Lemma 3.3 we have
\[
\alpha_{\ell}(\emph{\textbf{M}})\leq\|\gamma_{\ell}\|,
\ell=0,1,\ldots,k,
\]
where
\[
\gamma_{\ell}=(I-\max_{0\leq i\leq
k}\{\wedge^{-1}_{i}|C_{i}|\})^{-1}\wedge_{\ell}^{-1}.
\]
\end{theorem}
\textbf{Proof.} For the sake of convenience, we only prove the case when $\ell=0$, for $\ell=1,2,\ldots,k$, the proof is analogical. Now we show that
\begin{equation}\label{eq:313}
\alpha_{0}(\textbf{M})\leq\|(I-\max_{0\leq i\leq
k}\{\wedge^{-1}_{i}|C_{i}|\})^{-1}\wedge_{0}^{-1}\|.
\end{equation}

By the proof of Lemma 3.3, we have
\[
S_{\textbf{M}}^{-1}D_{0}=(V-U)^{-1}D_{0}=(I-V^{-1}U)^{-1}V^{-1}D_{0}
\]
and
\begin{align*}
|(I-V^{-1}U)^{-1}|\leq [I-\max_{0\leq i\leq k}\{\wedge^{-1}_{i}|C_{i}|\}]^{-1}.
\end{align*}
From Lemma 2.4 and
\[D_{1}\wedge_{1}+\ldots+D_{k}\wedge_{k}\geq(I-D_{0})
\min\{\wedge_{1},\wedge_{2},\ldots,\wedge_{k}\},
\]
we have
\begin{align*}
V^{-1}D_{0}=&[D_{0}\wedge_{0}+D_{1}\wedge_{1}+\ldots+D_{k}\wedge_{k}]^{-1}D_{0}\\
\leq&[D_{0}\wedge_{0}+(I-D_{0})\min\{\wedge_{1},\wedge_{2},\ldots,\wedge_{k}\}]^{-1}D_{0}\\
\leq&\wedge_{0}^{-1}.
\end{align*}
Hence, we have
\begin{align*}
\|S_{\textbf{M}}^{-1}D_{0}\|\leq&\||(I-V^{-1}U)^{-1}V^{-1}D_{0}|\|\\
\leq&\|[I-\max_{0\leq i\leq
k}\{\wedge^{-1}_{i}|C_{i}|\}]^{-1}\wedge^{-1}_{0}\|,
\end{align*}
which implies that the desired bound (\ref{eq:313}) holds. $\hfill{}
\Box$

For the VLCP ($\mathbf{M},\mathbf{q}$),  we have the following
corollary.

\begin{cor}
Let  all the diagonal entries of the matrices $M_{i}$ in
$\mathbf{M}=(I,M_{1},\ldots, M_{k})$ be positive. If
\begin{equation*}
\rho(\max_{1\leq i\leq k}\{\wedge^{-1}_{i}|C_{i}|\})<1,
\end{equation*}
then we have
\begin{description}
\item  $ (i)$ $\mathbf{M}=(I, M_{1},\ldots, M_{k})$ has the row $\mathcal{W}$-property.
\item  $ (ii)$ The following inequality holds:
\[
\tilde{\alpha}_{\ell}(\emph{\textbf{M}})\leq\tilde{\gamma}_{\ell},
\ell=1,2,\ldots,k,
\]
where
\[
\tilde{\gamma}_{\ell}= \|(I-\max_{1\leq i\leq
k}\{\wedge^{-1}_{i}|C_{i}|\})^{-1}\wedge_{\ell}^{-1}\|.
\]
\end{description}
\end{cor}

\begin{remark}
It is noted that for $k=1$  $M_{1}$ is an $H_+$-matrix in Corollary
3.1, Corollary 3.1(ii) directly reduces to Theorem 2.5 of
\cite{Chen07}.
\end{remark}
\subsubsection{Case 2}
Let $\mathbf{M}=(M_{0} ,M_{1},\ldots, M_{k})$. In the subsection, we
consider the case that $M_{i}$ in $\mathbf{M}$ is an sdd matrix for
$i=0,1,...,k$. First, we have

\begin{lemma}
Let $M_{i}$ in $\mathbf{M}$, $i=0,1,...,k$, be an sdd matrix with
the $s$-th diagonal entry having the same sign, $s=1,...,n$. Then
matrix $S_{\textbf{M}}$ for any $(D_{0},D_{1},\ldots,D_{k})\in
\mathcal{D}$ is also an sdd matrix.
\end{lemma}
\textbf{Proof.} Recall the notations given in Section 1. Then we set
\[
\langle M_{i}\rangle e=(r^{(i)}_{1}, r^{(i)}_{2},\ldots,
r^{(i)}_{n})^{T}.
\]
Since $M_{i}$ is an sdd matrix, we have $r_{i}=\min_{j}\{r^{(i)}_{j}\}>0,
j=1,2,\ldots,n$. By the simple computations we have $\langle S_{\textbf{M}}\rangle \geq \sum_{i=0}^{k}D_{i}\langle
M_{i}\rangle$. Hence
\[
\langle S_{\textbf{M}}\rangle e\geq\sum_{i=0}^{k}D_{i}\langle
M_{i}\rangle e\geq\sum_{i=0}^{k}D_{i}r_{i}e\geq
\min\{r_{i}\}e\sum_{i=0}^{k}D_{i}=\min\{r_{i}\}e>0,
\]
which shows that $S_{\textbf{M}}$ is an sdd matrix. $\hfill{} \Box$

\begin{theorem}
Under the assumption of Lemma 3.4 we have
\[
(\alpha_{\ell}(\emph{\textbf{M}}))_{\infty}\leq \delta_{\ell},
\ell=0,1,\ldots,k,
\]
where $(\alpha_{\ell}(\emph{\textbf{M}}))_{\infty}=\max_{\mathcal{D}}\|S_{\textbf{M}}^{-1}D_\ell\|_{\infty}$ and
\[
\delta_{\ell}=\frac{1}{\min_{i\in N}\{(\langle M_{\ell}\rangle
e)_{i}\}}.
\]
\end{theorem}
\textbf{Proof.} We only prove that
\begin{equation}\label{eq:314}
(\alpha_{0}(\textbf{M}))_{\infty}\leq\frac{1}{\min_{i\in
N}\{(\langle M_{0}\rangle e)_{i}\}}.
\end{equation}

By Lemma 3.4, $S_{\textbf{M}}$ is an sdd matrix. Let $D_{0}=\mbox
{diag}(d^{(0)}_1,...,d^{(0)}_n)$. It follows from Lemma 4 of
\cite{Li13} that
\begin{equation}\label{315}
\|S_{\textbf{M}}^{-1}D_0\|_{\infty}\leq \max_{i}\frac{d^{(0)}_i}{(\langle S_{\textbf{M}}\rangle e)_i}.
\end{equation}

Let $r=\min_{1\leq i\leq k;1\leq j\leq n}(\langle M_{i}\rangle
e)_j$. Since $M_i$ is an sdd matrix, we get $r>0$ and $\langle
M_{i}\rangle e\geq re$. By the proof of Lemma 3.4 we have

\[
\langle S_{\textbf{M}}\rangle e \geq
D_{0}\langle M_{0}\rangle e + \sum_{i=1}^{k}D_{i}\langle M_{i}\rangle e\geq D_{0}\langle M_{0}\rangle e +r\sum_{i=1}^{k}D_{i}e
=D_{0}\langle M_{0}\rangle e +r(I-D_{0})e.
\]
Hence,
\begin{equation}\label{316}
\frac{d^{(0)}_i}{(\langle S_{\textbf{M}}\rangle e)_i}\leq \frac{d^{(0)}_i}{d^{(0)}_i(\langle M_{0}\rangle e)_i +r(1-d^{(0)}_i)}.
\end{equation}
Then the desired bound \eqref{eq:314} follows from \eqref{315},
\eqref{316} and Lemma 2.4. This completes the proof of the theorem.
$\hfill{} \Box$

The following result is for the VLCP case.

\begin{cor}
Let $M_{i}$ in $\mathbf{M}=(I,M_{1},\ldots, M_{k})$ be an sdd matrix whit positive diagonals. Then
the following statements hold:
\begin{description}
\item  $ (i)$ $S^{v}_{\textbf{M}}$ for any
$(D_{0},D_{1},\ldots,D_{k})\in \mathcal{D}$  is an sdd matrix with positive diagonals, and $\mathbf{M}=(I, M_{1},\ldots, M_{k})$ has the row
$\mathcal{W}$-property.
\item  $ (ii)$  The following bound holds:
\[
(\tilde{\alpha}_{\ell}(\emph{\textbf{M}}))_{\infty}\leq
\tilde{\delta}_{\ell}, \ell=1,2,\ldots,k,
\]
where $(\tilde{\alpha}_{\ell}(\emph{\textbf{M}}))_{\infty}=\max_{\mathcal{D}}\|\tilde{S}_{\textbf{M}}^{-1}D_\ell\|_{\infty}$ and
\[
\tilde{\delta}_{\ell}=\frac{1}{\min_{i\in N}\{(\langle M_{\ell}
\rangle e)_{i}\}}.
\]
\end{description}
\end{cor}

\begin{remark}
The conditions in Lemmas 3.3  and 3.4 are not included each other,
e.g., taking the block matrix $\emph{\textbf{M}}=(M_{0}, M_{1})$ as
follows:
\[
M_{0}=\left[\begin{array}{ccc}
1&0\\
-1&1\\
\end{array}\right],\ M_{1}=\left[\begin{array}{ccc}
2&0\\
3&2\\
\end{array}\right].
\]

By the simple computation, we get
\[
\rho(\max\{\wedge^{-1}_{0}|C_{0}|,\wedge^{-1}_{1}|C_{1}|\})=0<1.
\]
Then $\emph{\textbf{M}}$ satisfies the condition in Lemma 3.3. But
$M_{0}$ and $M_{1}$ are not sdd, i.e.,  $\emph{\textbf{M}}$ does not
satisfy the condition in Lemma 3.4. Now we take
$\emph{\textbf{M}}=(M_{0}, M_{1})$ as follows:
\[
M_{0}=\left[\begin{array}{ccc}
2&1&0\\
1&2&0\\
0&1&2\\
\end{array}\right],\ M_{1}=\left[\begin{array}{ccc}
2&0&1\\
0&2&1\\
1&0&2\\
\end{array}\right].
\]
Then both $M_{0}$ and $M_{1}$ are sdd. However, one may get
\[
\rho(\max\{\wedge^{-1}_{0}|C_{0}|,\wedge^{-1}_{1}|C_{1}|\})=1,
\]
which shows that $\emph{\textbf{M}}$ does not satisfy the condition
in Lemma 3.3.
\end{remark}

\section{Relative perturbation bounds}

In this section, we discuss the relative perturbation bounds for
EVLCP ($\mathbf{M},\mathbf{q}$) and the VLCP
($\mathbf{M},\mathbf{q}$).

\begin{theorem}
Let $\emph{\textbf{M}}=(M_{0}, M_{1},\ldots, M_{k})$ have the row
$\mathcal{W}$-property, and let the perturbation $\Delta M_{i}\in
\mathbb{R}^{n\times n}$ and $\Delta q_{i}\in \mathbb{R}^{n}$ satisfy
$\|\Delta M_{i}\|\leq\epsilon_{i}\|M_{i}\|$ and $\|\Delta
q_{i}\|\leq\epsilon_{i}\|(-q_{i})_{+}\|$, respectively. If the
perturbation $\epsilon_{i}$ is small enough such that
$\eta:=\sum_{i=0}^{k}\epsilon_{i}\alpha_{i}(\textbf{\emph{M}})\|M_i\|<1$,
then we have:
\begin{description}
\item  $(i)$ $\tilde{\emph{\textbf{M}}}=(M_{0}+\Delta M_{0}, M_{1}+\Delta
M_{1},\ldots, M_{k}+\Delta M_{k})$ has the row
$\mathcal{W}$-property.
\item  $ (ii)$ Let $x$ and $y$ be the solution of EVLCP $(\mathbf{M},\mathbf{q})$ and EVLCP $(\tilde{\mathbf{M}},\tilde{\mathbf{q}})$, respectively. Then
\begin{eqnarray}\label{eq:41}
\frac{\|x-y\|}{\|x\|}\leq \frac{2\eta}{1-\eta}.
\end{eqnarray}
\end{description}
\end{theorem}
\textbf{Proof.} (i) follows immediately from Lemma 3.1.

(ii) By (\ref{36}) and (\ref{eq:39}) it is easy to show that
\[
\alpha_{i}(\tilde{\textbf{M}})\leq\frac{1}{1-\eta}\alpha_{i}(\textbf{M})
\]
and
\begin{eqnarray}\label{eq:42}
\|y-x\|\leq
\frac{1}{1-\eta}\bigg(\bigg(\sum_{i=0}^{k}\alpha_{i}(\textbf{M})\|\Delta
M_{i}\|\bigg)\|x\|+\sum_{i=0}^{k}\alpha_{i}(\textbf{M})\|\Delta
q_{i}\|\bigg).
\end{eqnarray}
From $M_{i}x+q_{i}\geq0$, we deduce that
\begin{equation*}
  (-q_{i})_{+}\leq(M_{i}x)_{+}\leq|M_{i}x|\leq|M_{i}||x|,
\end{equation*}
which together with the assumption gives
\begin{eqnarray}\label{l3}
\|\Delta
q_{i}\|\leq\epsilon_{i}\|(-q_{i})_{+}\|\leq\epsilon_{i}\|(M_{i}x)_{+}\|\leq\epsilon_{i}\|M_{i}\|\|x\|.
\end{eqnarray}
Combining (\ref{eq:42}) and \eqref{l3} together gives
\begin{align*}
\|y-x\|\leq&
\frac{1}{1-\eta}\bigg(\sum_{i=0}^{k}\alpha_{i}(\textbf{\emph{M}})\|\Delta
M_{i}\|+\sum_{i=0}^{k}\alpha_{i}(\textbf{M})\epsilon_{i}\|M_{i}\|\bigg)\|x\|\\
\leq&
\frac{2}{1-\eta}\bigg(\sum_{i=0}^{k}\epsilon_{i}\alpha_{i}(\textbf{\emph{M}})\|M_{i}\|\bigg)\|x\|,
\end{align*}
from which one may deduce the desired bound \eqref{eq:41}. $\hfill{}
\Box$

It is noticed that  VLCP is a special case of EVLCP. The relative
perturbation bound for VLCP $(\mathbf{M},\mathbf{q})$ can be deduced
from Theorem 4.1. Here we omit it.

\begin{remark}
If we take $k=1,~q_0=0,~M_0=I, \mbox{and}~M_{1}$ is a $P$-matrix in Theorem 4.1,
then Theorem 3.1 of \cite{Chen07} can be derived from Theorem 4.1.
\end{remark}

In the following, we consider the special cases where $M_{i},
i=0,1,...,k$, has positive diagonals or an sdd matrix as in the
subsections 3.2.1 and 3.2.2, see Corollary 4.1. Its proof is similar
to Theorem 4.1. We omit it.

\begin{cor}
Under the same assumption as in Theorem 4.1, if $M_{i}$ satisfies
the same assumption as in Lemma 3.3 or Lemma 3.4 for any
$i=0,1,...,k$, respectively, then
$\tilde{{\textbf{\emph{M}}}}=(M_{0}+\Delta M_{0}, M_{1}+\Delta
M_{1},\ldots, M_{k}+\Delta M_{k})$ has the row
$\mathcal{W}$-property and
\begin{eqnarray}\label{b1}
\frac{\|x-y\|}{\|x\|}\leq \frac{2\eta_{\xi}}{1-\eta_{\xi}},
\end{eqnarray}
provided that the perturbation $\epsilon_{i}$ is small enough such
that $\sum_{i=0}^{k}\epsilon_{i}\|\xi_i\|\|M_i\|= \eta_{\xi}< 1$,
where $\xi=\gamma~\mbox{or}~\delta$ is given by Theorem 3.4 or
Theorem 3.5, respectively.
\end{cor}

It is noted that $\|\delta_i\|=\delta_i$ in the bound (\ref{b1})
(see Theorem 3.5). Next, we consider a special norm for the relative
perturbation bound. The proof is tedious, we omit it.

\begin{theorem} Let $\emph{\textbf{M}}=(M_{0}, M_{1},\ldots, M_{k})$ have the row
$\mathcal{W}$-property,  the perturbation $\Delta M_{i}\in
\mathbb{R}^{n\times n}$ and $\Delta q_{i}\in \mathbb{R}^{n}$ satisfy
$|\Delta M_{i}|\leq\epsilon_{i}|M_{i}|$ and $|\Delta
q_{i}|\leq\epsilon_{i}|(-q_{i})_{+}|$, respectively. Let
$M_{i}=\wedge_{i}-C_{i}$ be an $H_{+}$-matrix and $\rho(\max_{0\leq
i\leq k}\{\wedge^{-1}_{i}|C_{i}|\})<1$,  where $\wedge_{i}$ is the
diagonal part of $M_{i}$, $i=0,1,\ldots,k$. If the perturbation
$\epsilon_{i}$ is small enough such that
\[
\sum_{i=0}^{k}\epsilon_{i}\|\gamma_{i}\|_{\infty}\|M_{i}\|_{\infty}=\hat{\eta}<1,
\]
then the following statements hold:
\begin{description}
\item  $ (i)$ $\tilde{\emph{\textbf{M}}}=(M_{0}+\Delta
M_{0}, M_{1}+\Delta M_{1},\ldots, M_{k}+\Delta M_{k})$ has the row
$\mathcal{W}$-property.
\item  $ (ii)$ Let $x$ and $y$ be the solution of EVLCP ($\mathbf{M},\mathbf{q}$) and EVLCP $(\tilde{\mathbf{M}},\tilde{\mathbf{q}})$, respectively. Then
\begin{eqnarray}\label{eq:43}
\frac{\|x-y\|_{\infty}}{\|x\|_{\infty}}\leq\frac{2(\epsilon_{0}\|\gamma_{0}|M_{0}|\|_{\infty}+\ldots+\epsilon_{k}\|\gamma_{k}|M_{k}|\|_{\infty})}{1-\hat{\eta}}.
\end{eqnarray}
\end{description}
\end{theorem}

\begin{remark}
Taking $k=1,~M_0=I$ and $q_0=0$, the bound \eqref{eq:43}  reduces to the corresponding one in Theorem 3.3 of \cite{Chen07}.
\end{remark}

\section{Numerical examples}
In this section, some numerical examples are given to show the
feasibility of the relative perturbation bound. For the sake of
convenience, we only use the infinity norm in all numerical
experiments. Let $x$ and $y$, respectively, be the solution of EVLCP
(\ref{p1}) and the perturbed EVLCP (1.4), which can be obtained by
directly using  Lemke's complementarity pivoting method
\cite{Cottle70} for the following examples, and $r$ be the real
relative error given by
\[
r=\frac{\|x-y\|_{\infty}}{\|x\|_{\infty}}.
\]

The perturbation for EVLCP (\ref{p1}) can be set as:
\[
\Delta
M_{i}=\frac{\epsilon\|M_{i}\|_{\infty}}{\|S_{i}\|_{\infty}}S_{i},
\Delta
q_{i}=\frac{\epsilon\|q_{i}\|_{\infty}}{\|t_{i}\|_{\infty}}t_{i},
\]
with $S_{i}$ and $t_{i}$,  respectively, being an arbitrary random
matrix and vector. In this case we have  $$\|\Delta M_{i}\|_{\infty}\leq
\epsilon \|M_{i}\|_{\infty},  ~\|\Delta
q_{i}\|_{\infty}\leq\epsilon\|q_{i}\|_{\infty}.$$ 
All the computations are done in Matlab R2021b.

\textbf{Example 5.1} (\cite{Qi99,Ebiefung95}) Let $k=1$ in
(\ref{p1}), $\textbf{M}=(M_{0}, M_{1})$ and $\mathbf{q}=(q_{1},
q_{2})$, respectively, be form
\[
M_{0}=\left[\begin{array}{ccc}
3&-2\\
-4&5\\
\end{array}\right],\ M_{1}=\left[\begin{array}{ccc}
4&-4\\
-1&2\\
\end{array}\right], \ \mbox{and}\ q_{0}=q_{1}=(-1,-1)^{T}.
\]
It is easy to check that $\textbf{M}$ has the row
$\mathcal{W}$-property. This implies that the corresponding EVLCP
($\mathbf{M},\mathbf{q}$) has a unique solution. In fact, its unique
solution $x=(2.25,2)^{T}$.

\begin{table}[!htb] \centering
\begin{tabular}
{p{50pt}p{100pt}p{55pt}p{50pt}p{50pt}p{50pt}} \hline
$\epsilon$&$y$&$r$&$\bar{\tau}_{\gamma}$&$\tau$& $\nu$\\\hline
0.01& $(1.9277, 1.7120)^{T}$&0.1444& 1.2359 &   2.3478 & 1.5145  \\
0.001&$(2.2132,1.9671)^{T}$&0.0155  & 0.0518&0.1142  &0.0736  \\
0.0001&$(2.2463,1.9967)^{T}$&0.0017 &0.0057 & 0.0109&0.0070  \\
\hline
\end{tabular}
\\ \caption{Relative perturbation bounds of Example 5.1.}
\end{table}

By the simple computations,  we obtain that
\[
\rho(\max\{\wedge^{-1}_{0}|C_{0}|,\wedge^{-1}_{1}|C_{1}|\})=0.8944<1.
\]
This shows that the condition in Theorem 4.2 and Corollary 4.1 with
Lemma 3.3 are satisfied. Since $M_{1}$ is not a strictly row
diagonally dominant matrix, the conditions of Corollary 4.1 with
Lemma 3.4 are not satisfied. Based on (\ref{eq:42}), Theorem 4.2 and
Corollary 4.1 with Lemma 3.3, we set
$\eta_{\gamma}=\epsilon\sum_{i=0}^{1}\|\gamma_i\|_{\infty}\|M_i\|_{\infty}$,
\[
\bar{\tau}_{\gamma}=\frac{1}{1-\eta_{\gamma}}\bigg(\sum_{i=0}^{1}\|\gamma_{i}\|_{\infty}\|\Delta
M_{i}\|_{\infty}+\frac{\sum_{i=0}^{1}\|\gamma_{i}\|_{\infty}\|\Delta
q_{i}\|_{\infty}}{\|x\|_{\infty}}\bigg),
\]
and
\[
\tau=\frac{2\eta_{\gamma}}{1-\eta_{\gamma}},
\nu=\frac{2\epsilon(\|\gamma_{0}|M_{0}|\|_{\infty}+\|\gamma_{1}|M_{1}|\|_{\infty})}{1-\eta_{\gamma}}.
\]
In our computations, we choose some values of $\epsilon$ such that
$\eta_{\gamma}<1$, see Table 1.

We report the numerical result for three relative perturbation
bounds in Table 1, from which we find that
$\bar{\tau}_{\gamma}<\tau$ and $\bar{\tau}_{\gamma}<\nu$, and also
illustrates that the proposed bounds are very close to the real
relative value when the perturbation is very small. This show that
the relative perturbation bounds given by Corollary 4.1 and Theorem
4.2 are feasible and effective under some suitable condition. The
following example is given by \cite{Zhang09}.

\textbf{Example 5.2} Consider the EVLCP ($\mathbf{M},\mathbf{q}$),
in which is given $\textbf{M}=(M_{0}, M_{1})$,  where
\[
M_{0}=\left[\begin{array}{ccc}
1&3/4&0\\
3/4&1&0\\
0&3/4&1\\
\end{array}\right],\ M_{1}=\left[\begin{array}{ccc}
1&0&3/4\\
0&1&3/4\\
3/4&0&1\\
\end{array}\right],
\]
Clearly, $M_{0}$ and $M_{1}$ are two strictly row diagonally
dominant matrices. Hence, EVLCP  has a unique solution for any
$\mathbf{q}$ (see Lemma 3.4). Let
$q_{1}=q_{2}=(-1.75,-1.75,-1.75)^{T}$. Then it is easy to see that
$x=(1,1,1)^{T}$ is its unique solution.

\begin{table}[!htb] \centering
\begin{tabular}
{p{50pt}p{150pt}p{55pt}p{55pt}p{55pt}} \hline
$\epsilon$&$y$&$r$&$\bar{\tau}_{\delta}$&$\upsilon$\\\hline
0.01& $( 1.0510,  1.0545,  0.9615)^{T}$& 0.0545& 0.3256&  0.3256 \\
0.001&$(  1.0078,  1.0076,  0.9930)^{T}$&  0.0078& 0.0284& 0.0284  \\
0.0001&$( 1.0007, 1.0007, 0.9994)^{T}$& 7.0368e-04&0.0028 &0.0028\\
\hline
\end{tabular}
\\ \caption{Relative perturbation  bounds of Example 5.2.}
\end{table}

Clearly, the condition of  Corollary 4.1 with Lemma 3.4 is
satisfied. However, $\textbf{M}$ does not satisfy the condition of
Theorem 4.2 and Corollary 4.1 with Lemma 3.3 because
\[
\rho(\max\{\wedge^{-1}_{0}|C_{0}|,\wedge^{-1}_{1}|C_{1}|\})=1.5>1.
\]
Based on (\ref{eq:42}) and  Corollary 4.1 with Lemma 3.4,  we set
$\eta_{\delta}=\epsilon\sum_{i=0}^{1}\delta_i\|M_i\|_{\infty}$,
\[
\bar{\tau}_{\delta}=\frac{1}{1-\eta_{\delta}}\bigg(\sum_{i=0}^{1}\delta_i\|\Delta
M_{i}\|_{\infty}+\frac{\sum_{i=0}^{1}\delta_i\|\Delta
q_{i}\|_{\infty}}{\|x\|_{\infty}}\bigg),\ \mbox{and}\
\upsilon=\frac{2\eta_{\delta}}{1-\eta_{\delta}}.
\]
We take some values of $\epsilon$ such that $\eta_{\delta}<1$.

The numerical bounds are reported in Table 2, which show that
$\bar{\tau}_{\delta}=\upsilon$. It also illustrates that the
numerical bounds in Example 5.2 show the same perturbing behavior as
in Example 5.1 although conditions are different.

\vspace{0.3cm} Next, we give a example from the discretization of
Hamilton-Jacobi-Bellman (HJB) equation, in which the conditions in
both Corollary 4.1 and Theorem 4.3 hold. \vspace{0.3cm}

\textbf{Example 5.3} Consider the following EVLCP
($\mathbf{M},\mathbf{q}$)
\[
\min\{M_{0}x+q_{0}, M_{1}x+q_{1}\}=0,
\]
where $M_{i}$ and $q_{i} \ (i=0,1)$ comes from the discretization of
Hamilton-Jacobi-Bellman (HJB) equation
\begin{equation*}
\left\{ \begin{aligned} &\max_{0\leq i\leq 1}\{L_{i}+f_{i}\}=0 \ \rm{in} \ \Gamma, \\
&u=0 \ \rm{on}\  \partial\Gamma,
\end{aligned} \right.
\end{equation*}
with $\Gamma=\{(x,y)|0<x<2,0<y<1\}$,
\begin{equation*}
\left\{ \begin{aligned}
&L_{0}=0.002u_{xx}+0.001u_{yy}-20u,f_{0}=1, \\
&L_{1}=0.001u_{xx}+0.001u_{yy}-10u, f_{1}=1,
\end{aligned} \right.
\end{equation*}
see \cite{Bensoussan82} for more details. Here, by making use of the
central difference scheme to discretize the above HJB equation,
Example 5.3 can be obtained and $q_{0}=q_{1}=-e$.

Here, $M_{0}$ and $M_{1}$ obtained from the above HJB equation are
two strictly row diagonally dominant matrices. What's more,
\[
\rho(\max\{\wedge^{-1}_{0}|C_{0}|,\wedge^{-1}_{1}|C_{1}|\})\leq\|W\|_{\infty}<1,
\]
where $W=\max\{\wedge^{-1}_{0}|C_{0}|,\wedge^{-1}_{1}|C_{1}|\}$.
This means that the conditions in Corollary 4.1 are satisfied, so do
the conditions in Theorem 4.2. This implies that Example 5.3 has a
unique solution because the block matrix $\textbf{M}$ has the row
$\mathcal{W}$-property.

Tables 3-6 list some relative perturbation bounds for Example 5.3
with the different dimension and $\epsilon$, in which $r$, $\tau$,
$\nu$, $\upsilon$, $\eta_{\gamma}$ and $\eta_{\delta}$ are the above
defined in  Example 5.1 and Example 5.2. Again, we chose some values
of $\epsilon$ such that $\eta_{\gamma}<1$ and $\eta_{\delta}<1$.
(a), (b), (c) and (d) in Figure 1 are in line with Table 3, Table 4,
Table 5 and Table 6, respectively. In Figure 1, `RPB' denotes the
value of the relative perturbation bound and `$n$' denotes the order
of the system matrix.

\begin{table}[!htb] \centering
\begin{tabular}
{p{50pt}p{55pt}p{55pt}p{55pt}p{55pt}} \hline
$\epsilon$&$r$&$\tau$&$\upsilon$ &$\nu$\\\hline
0.01&0.0175  &0.0500 &0.0454 &0.0499 \\
0.015&0.0273 &0.0760 &0.0688 &0.0759 \\
0.02&0.0356  &0.1026 &0.0929 &0.1025 \\
0.025&0.0433  &0.1300 &0.1174 &0.1297 \\
0.03&0.0547  &0.1580 &0.1426 &0.1577 \\
\hline
\end{tabular}
\\ \caption{Relative perturbation bounds of Example 5.3 with $n=16$.}
\end{table}

\begin{table}[!htb] \centering
\begin{tabular}
{p{50pt}p{55pt}p{55pt}p{55pt}p{55pt}} \hline
$\epsilon$&$r$&$\tau$&$\upsilon$ &$\nu$\\\hline
0.01& 0.0191 & 0.0551& 0.0476&0.0551 \\
0.015&0.0281  &0.0838 &0.0722 &0.0838 \\
0.02&0.0373  &0.1133 &0.0975 &0.1133 \\
0.025&0.0460  & 0.1437&0.1234 &0.1436 \\
0.03& 0.0564 &0.1749 &0.1499 &0.1749 \\
\hline
\end{tabular}
\\ \caption{Relative perturbation bounds of Example 5.3 with $n=36$.}
\end{table}

\begin{table}[!htb] \centering
\begin{tabular}
{p{50pt}p{55pt}p{55pt}p{55pt}p{55pt}} \hline
$\epsilon$&$r$&$\tau$&$\upsilon$ &$\nu$\\\hline
0.01&0.0193  &0.0609 &0.0500 & 0.0609\\
0.015&0.0292  &0.0927 &0.0760 & 0.0927\\
0.02& 0.0375 &0.1255 &0.1026 &0.1255 \\
0.025& 0.0472 &0.1594 &0.1299 &0.1594 \\
0.03& 0.0576 &0.1944 &0.1580 &0.1944 \\
\hline
\end{tabular}
\\ \caption{Relative perturbation bounds of Example 5.3 with $n=64$.}
\end{table}

\begin{table}[!htb] \centering
\begin{tabular}
{p{50pt}p{55pt}p{55pt}p{55pt}p{55pt}} \hline
$\epsilon$&$r$&$\tau$&$\upsilon$ &$\nu$\\\hline
0.01&0.0199  &0.0676 &0.0527 &0.0676 \\
0.015&0.0295  &0.1032 &0.0801 &0.1032  \\
0.02&0.0390  &0.1400 &0.1083 &0.1400 \\
0.025&0.0480  &0.1782 &0.1372 &0.1782 \\
0.03&0.0585  &0.2177 &0.1670 &0.2177 \\
\hline
\end{tabular}
\\ \caption{Relative perturbation bounds of Example 5.3 with $n=100$.}
\end{table}

\begin{figure}
\setcaptionwidth{5in} \subfigure[$n=16$]{
\begin{minipage}[b]{0.5\textwidth}
\centering
\includegraphics[width=2.5in]{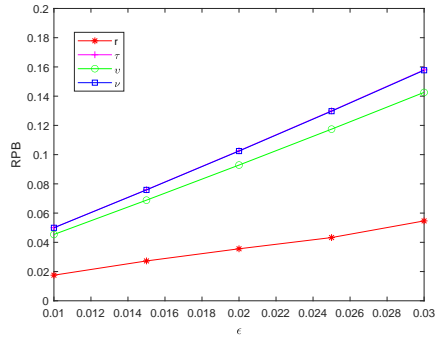}
\end{minipage}}%
 \subfigure[$n=36$]{
\begin{minipage}[b]{0.5\textwidth}
\centering
\includegraphics[width=2.5in]{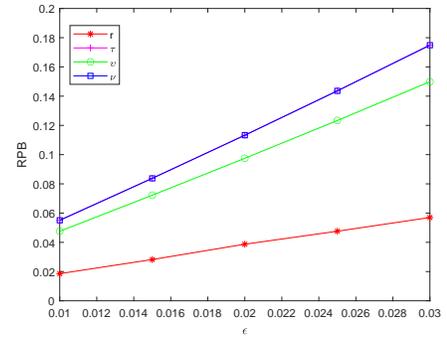}
\end{minipage}}
\subfigure[$n=64$]{
\begin{minipage}[b]{0.5\textwidth}
\centering
\includegraphics[width=2.5in]{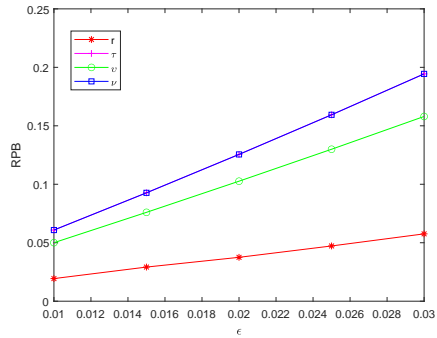}
\end{minipage}}
\subfigure[$n=100$]{
\begin{minipage}[b]{0.5\textwidth}
\centering
\includegraphics[width=2.5in]{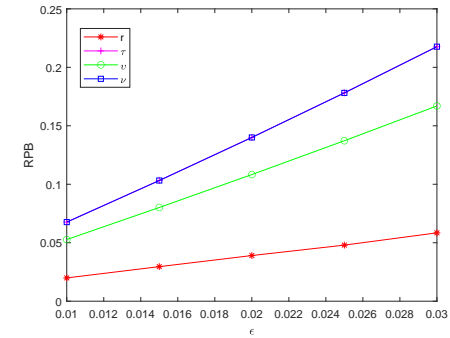}
\end{minipage}}
\caption{The value and number of solutions of Example 5.3.}
\end{figure}

Similar to what happens in Example 5.1 and Example 5.2, from Tables
3-6, we can draw the same conclusion. In other word, for the same
dimension, with $\epsilon$ decreasing, $r$, $\upsilon$, $\tau$ and
$\nu$ are decreasing, also see Figure 1. The reason is the same as
Example 5.1 and Example 5.2. In addition, we find that for the same
$\epsilon$, with the dimension increasing, $\upsilon$, $\tau$ and
$\nu$ are increasing (the reason is similar to the same dimension
with $\epsilon$ decreasing), and the values of $\tau$ and $\nu$ are
fairly close, not much different in size.

No matter what, from the above numerical results in Tables 3-6, we
still verify that under some suitable condition, Corollary 4.1 and
Theorem 4.2 indeed provide some valid relative perturbation bounds.

\section{Conclusion}
In this paper, we  discuss the perturbation analysis  of the EVLCP
($\mathbf{M},\mathbf{q}$). By making use of  a general equivalent
form of the minimum function, under the assumption of row $\mathcal{W}$-property,
some perturbation bounds for the EVLCP ($\mathbf{M},\mathbf{q}$)
are presented, which cover some existing results in
\cite{Chen07}. Particularly, for all diagonal elements of the
matrices $M_{i}$ in $\mathbf{M}$ being positive and all the matrices
$M_{i}$ in $\mathbf{M}$ being a strictly row diagonally dominant
matrix, some computable perturbation bounds are provided as well. Some numerical examples are given
to show the proposed bounds.

\section*{Acknowledgements} The authors would like to thank two anonymous
referees for providing helpful suggestions, which greatly improved
the paper.

\end{document}